\documentclass[12pt]{article}

\usepackage{amsmath,amssymb,amsthm,fullpage}
\usepackage{multirow}

\newtheorem{theorem}{Theorem}
\newtheorem{corollary}[theorem]{Corollary}

\title{On Avoiding Sufficiently Long Abelian Squares}
 
\author{Elyot Grant}
 
\begin{document}

\maketitle

\begin{abstract}
A finite word $w$ is an \emph{abelian square} if $w = xx^\prime$ with $x^\prime$ a permutation of $x$.  In 1972, Entringer, Jackson, and Schatz proved that every binary word of length $k^2 + 6k$ contains an abelian square of length $\geq 2k$.  We use Cartesian lattice paths to characterize abelian squares in binary sequences, and construct a binary word of length $q(q+1)$ avoiding abelian squares of length $\geq 2\sqrt{2q(q+1)}$ or greater.  We thus prove that the length of the longest binary word avoiding abelian squares of length $2k$ is $\Theta(k^2)$.
\end{abstract}

\section{Introduction}
Let $\Sigma$ be a finite alphabet.  A word $w \in \Sigma^*$ is an \emph{abelian square of order $k$} if $w = xx^\prime$ with $|x| = |x^\prime| = k$ and $x^\prime$ a permutation of $x$.  In 1972, Entringer, Jackson, and Schatz proved that all infinite binary sequences contain arbitrarily large abelian squares \cite{EJS}.  In particular, they showed that all binary words $w \in \{0,1\}^*$ of length $k^2 + 6k$ contain an abelian square of order $k$ or greater.  In this paper, we examine $\ell(k)$, the length of the longest binary word avoiding abelian squares $xx^\prime$ with $|x| \geq k$.  

Precise values of $\ell(k)$ have been computed for $1 \leq k \leq 10$ by Jeffrey Shallit and Narad Rampersad via a brute force search.  The results are given in Section~\ref{secvals}.

The bound $\ell(k) < k^2 + 6k$ given by Entringer, Jackson, and Schatz is not the best possible upper bound, but an improved upper bound remains unknown.  A simple lower bound $\ell(k) \geq 8k - 6$ can be obtained by observing that the string $0^{2k-2}1^{2k-1}0^{2k-1}1^{2k-2}$ contains no abelian squares of order $k$ or greater.  This lower bound is tight for $2 \leq k \leq 7$, but is suboptimal for $k \geq 8$.

In this paper, we give a quadratic lower bound for $\ell(k)$, proving that $\ell(k)$ is $\Theta(k^2)$.  Moreover, we provide an intuitive geometric characterization of abelian squares in a binary word by treating each character of a string as a step of a lattice path in the Cartesian plane.  We use this geometric notion to construct, for all $q$, a word of length $q(q+1)$ containing no abelian squares of order $\geq \sqrt{2q(q+1)}$.

Many thanks go to Jeffrey Shallit for suggesting this as a problem to study as part of \emph{CS~860: Patterns in Strings: Existence, Avoidability, Enumeration}, a course he developed and taught at the University of Waterloo.

\section{Values of $\ell(k)$ for $1 \leq k \leq 10$} \label{secvals}
Jeffrey Shallit and Narad Rampersad have provided the values of $\ell(k)$ for $1 \leq k \leq 10$.  We give them here, alongside the lexicographically least word of length $\ell(k)$ containing no abelian squares of order $k$ or greater:

\begin{center}
  \begin{tabular}{ | c | c | p{12.8cm} | }
    \hline
    $k$ & $\ell(k)$ & \\ \hline
    1 & 3 & 010 \\ \hline
    2 & 10 & 0011100011 \\ \hline
    3 & 18 & 000011111000001111 \\ \hline
    4 & 26 & 00000011111110000000111111 \\ \hline
    5 & 34 & 0000000011111111100000000011111111 \\ \hline
    6 & 42 & 000000000011111111111000000000001111111111 \\ \hline
    7 & 50 & 00000000000001000001100001111001111101111111111111 \\ \hline
    8 & 62 & 00000000000000010000100100011001100111011011110111111111111111 \\ \hline
    9 & 76 & 00000000000000000100000001100100001110100011110110011111110111
             11111111111111 \\ \hline
    10 & 90 & 00000000000000000001000000100100000110101000011110101001111101
              1011111101111111111111111111 \\ \hline
    11 & $\geq 106$ &  \\ \hline
    12 & $\geq 124$ &  \\ \hline
    13 & $\geq 139$ &  \\ \hline
  \end{tabular}
\end{center}
 
\section{Main Result}
Given a word $w[1..t] \in \{0,1\}^*$, let $S_i = \sum_{j=1}^{i}{w[j]}$ be the nondecreasing sequence of \emph{prefix sums} of $w$.  By plotting the ordered pairs $(i,S_i)$ for $0 \leq i \leq t$, we obtain a representation of $w$ as a path across the Cartesian lattice, stepping east when $w$ contains a zero, and northeast when $w$ contains a 1.  An example for the string $100110001$ is shown in Figure~1.

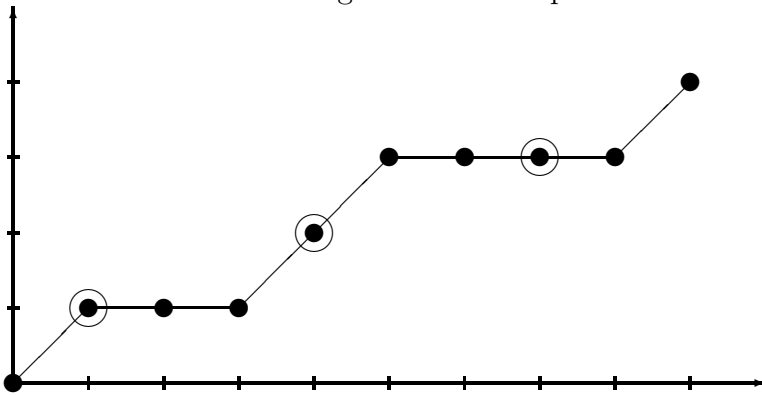
\begin{figure}
\caption{Lattice path for 100110001}
\setlength{\unitlength}{10mm}
\begin{picture}(11,6)(-1,-1)

%lines
\put(0,0){\line(1,1){1}}
\put(1,1){\line(1,0){1}}
\put(2,1){\line(1,0){1}}
\put(3,1){\line(1,1){1}}
\put(4,2){\line(1,1){1}}
\put(5,3){\line(1,0){1}}
\put(6,3){\line(1,0){1}}
\put(7,3){\line(1,0){1}}
\put(8,3){\line(1,1){1}}

%points
\put(0,0){\circle*{0.25}}
\put(1,1){\circle*{0.25}}
\put(2,1){\circle*{0.25}}
\put(3,1){\circle*{0.25}}
\put(4,2){\circle*{0.25}}
\put(5,3){\circle*{0.25}}
\put(6,3){\circle*{0.25}}
\put(7,3){\circle*{0.25}}
\put(8,3){\circle*{0.25}}
\put(9,4){\circle*{0.25}}

%circled points
\put(1,1){\circle{0.5}}
\put(4,2){\circle{0.5}}
\put(7,3){\circle{0.5}}

%axes
\linethickness{1pt}
\put(0,0){\vector(1,0){10}}
\put(0,0){\vector(0,1){5}}

%tickmarks
\put(1,-0.08){\line(0,1){0.16}}
\put(2,-0.08){\line(0,1){0.16}}
\put(3,-0.08){\line(0,1){0.16}}
\put(4,-0.08){\line(0,1){0.16}}
\put(5,-0.08){\line(0,1){0.16}}
\put(6,-0.08){\line(0,1){0.16}}
\put(7,-0.08){\line(0,1){0.16}}
\put(8,-0.08){\line(0,1){0.16}}
\put(9,-0.08){\line(0,1){0.16}}
\put(-0.08,1){\line(1,0){0.16}}
\put(-0.08,2){\line(1,0){0.16}}
\put(-0.08,3){\line(1,0){0.16}}
\put(-0.08,4){\line(1,0){0.16}}
\end{picture}
\end{figure}

We note that the number of ones in $w[m..n]$ is $S_n - S_{m-1}$.  Consequently, $w[i+1..i+2r]$ is an abelian square iff $S_{i+r} - S_i = S_{i+2r} - S_{i+r}$, which occurs precisely when $(i,S_i)$, $(i+r,S_{i+r})$, and $(i+2r,S_{i+2r})$ are three equally spaced collinear points in our lattice path.  In Figure~1, the three circled points indicate the presence of the subword $001100$, an abelian square.

Next, we give our construction of a word of length $q(q+1)$ containing no abelian squares of order $\geq \sqrt{2q(q+1)}$.  We design our word $w$ so that its lattice path approximates a quadratic function; this ensures that three equally spaced points along the path can be collinear only if they are sufficiently close together.  For $0 \leq i \leq q(q+1)$, define

\[
a_i = \left\lfloor \frac{i^2}{2q(q+1)} \right\rfloor \text{.}
\]

We note that if $i \leq q(q+1)$, then $i^2 - (i-1)^2 = 2i-1 < 2q(q+1)$, and hence $a_i - a_{i-1} \in \left\{0,1\right\}$ for all $1 \leq i \leq q(q+1)$. We can thus define a binary word $w = w[1..q(q+1)]$ by $w[i] = a_i - a_{i-1}$.  We will show the following:

\begin{theorem}\label{thetheorem}
$w$ contains no abelian squares $xx^\prime$ with $|x| \geq \sqrt{2q(q+1)}$.
\end{theorem}

Our theorem implies that if $q$ is an integer with $2q(q+1) \leq k^2$, then there exists a binary word of length $q(q+1)$ containing no abelian squares of order $k$.  For a given $k$, the shortest such $q$ is $\left\lfloor \frac{\sqrt{1+2k^2}-1}{2} \right\rfloor$.  Consequently, we may conclude the following:

\begin{corollary}
$\ell(k) \geq \left(\left\lfloor \displaystyle\frac{\sqrt{1+2k^2}-1}{2} \right\rfloor\right) \left(\left\lfloor \displaystyle\frac{\sqrt{1+2k^2}-1}{2} \right\rfloor + 1\right) > \displaystyle\frac{k^2}{2} - \sqrt{2k}$.
\end{corollary}

\section{Proof of Theorem \ref{thetheorem}}
Suppose $w$ contains an abelian square $xx^\prime$ with $|x| = r$.  Then there exist two adjacent blocks $w[i+1..i+r]$ and $w[i+r+1..i+2r]$ such that $|w[i+1..i+r]|_1 = |w[i+r+1..i+2r]|_1$.  This implies that $a_{i+r} - a_i = a_{i+2r} - a_{i+r}$.  We eliminate the floor function to bound the various $a_i$ values above and below in the following manner:
\begin{align*}
\frac{i^2}{2q(q+1)} - 1 &< a_i \\
a_{i+r} &\leq \frac{(i+r)^2}{2q(q+1)} \\
\frac{(i+2r)^2}{2q(q+1)} - 1 &< a_{i+2r} \\
\end{align*}
Taking a linear combination of the above inequalities, we obtain
\[
\frac{i^2}{2q(q+1)} - 1 + 2a_{i+r} + \frac{(i+2r)^2}{2q(q+1)} - 1 < a_i + 2\frac{(i+r)^2}{2q(q+1)} + a_{i+2r}
\]
and we may cancel the $a_i$ terms since $a_{i+r} - a_i = a_{i+2r} - a_{i+r}$.  We simplify what remains to obtain our result:
\begin{align*}
\frac{i^2}{2q(q+1)} + \frac{(i+2r)^2}{2q(q+1)} - 2 &< 2\frac{(i+r)^2}{2q(q+1)}\\
i^2 + (i+2r)^2 - 4q(q+1) &< 2(i+r)^2\\
r^2 &< 2q(q+1)\\
\end{align*}

\section{Additional Remarks}
One might suggest that we could improve our lower bound slightly by computing more $a_i$ values and extending $w$ to a longer string.  Indeed, we can take $a_i = \left\lfloor \frac{i^2}{2q(q+1)} \right\rfloor$ for all $i$ until we reach an $n$ such that $a_{n+1} - a_n > 1$.  Unfortunately, it turns out that this doesn't help us much.  Taking $p = {q(q+1) + \left\lceil \sqrt{2q(q+1)} \right\rceil}$, we see that
\begin{align*}
a_p - a_{q(q+1)} &= \left\lfloor \frac{(q(q+1) + \left\lceil \sqrt{2q(q+1)} \right\rceil)^2}{2q(q+1)} \right\rfloor -
                    \left\lfloor \frac{(q(q+1))^2}{2q(q+1)} \right\rfloor \\
                 &= \left\lceil \sqrt{2q(q+1)} \right\rceil + \frac{(\left\lceil \sqrt{2q(q+1)} \right\rceil)^2}{2q(q+1)} \\
                 &\geq p - q(q+1) + 1 \text{.}
\end{align*}
Consequently, there must be some $n$ with $q(q+1) \leq n < p$ such that $a_{n+1} - a_n > 1$.  Thus we can extend $w$ for at most another $\sqrt{2q(q+1)}$ symbols.

\end{document}